\documentclass[11pt,letterpaper]{article}
\usepackage{fullpage,tikz}
\usepackage{amsmath,amssymb}

\newcommand{\eg}{{\it e.g.}}
\newcommand{\ie}{{\it i.e.}}
\newcommand{\complex}{\mathbb{C}}
\newcommand{\reals}{\mathbb{R}}

\newcommand{\Tr}{\mathop{\bf tr}}
\newcommand{\diag}{\mathop{\bf diag}}
\newcommand{\Lflow}{\mathcal{F}}
\newcommand{\etal}{\emph{et al.}}
\newcommand{\V}{v}
\newcommand{\I}{i}
\newcommand{\iu}{\jmath}
\newcommand{\Pd}{P^{\mathrm{d}}}
\newcommand{\Tb}{\widetilde T}
\newcommand{\Yk}{Y_k}
\newcommand{\Ybk}{\widetilde Y_k}
\newcommand{\Pg}{p}
\newcommand{\Qd}{Q^{\mathrm{d}}}
\newcommand{\Qg}{q}
\newcommand{\Sd}{S^{\mathrm{d}}}
\newcommand{\Sg}{s}

\title{Reduced-Complexity Semidefinite Relaxations of\\Optimal Power Flow Problems}
\author{ Martin S.\ Andersen\thanks{Martin S. Andersen is with the Department of Applied
    Mathematics and Computer Science, Technical University of Denmark,
    2800 Lyngby, Denmark (e-mail: mskan@dtu.dk). Part of this research
    was conducted while he was a postdoc at Link\"oping University.
    Research supported by ERC Grant ERC-2011-AdG-20110209.} 
\and Anders Hansson\thanks{Anders Hansson is with the Division of Automatic Control,
    Department of Electrical Engineering, Link\"oping University,
    SE-581 83 Link\"oping, Sweden (e-mail:
    anders.g.hansson@liu.se). Research supported by the Swedish
    Department of Education within the ELLIIT project.}%
\and Lieven Vandenberghe\thanks{Lieven Vandenberghe is with the Electrical Engineering
    Department, University of California, Los Angeles, 66-147L
    Engineering IV, CA 90095-1594 (e-mail:
    vandenbe@ee.ucla.edu). Research supported by NSF Grants
    DMS-1115963 and ECCS-1128817.}}
\date{}

\begin{document}
\maketitle
\begin{abstract}
We propose a new method for generating semidefinite relaxations of
  optimal power flow problems. The method is based on chordal
  conversion techniques: by dropping some equality constraints
  in the conversion, we obtain semidefinite relaxations that
  are computationally cheaper, but potentially weaker, than the standard
  semidefinite relaxation. Our numerical results show that the new
  relaxations often produce the same results as the standard
  semidefinite relaxation, but at a lower computational cost.
\end{abstract}


\section{Introduction}
The general AC optimal power flow (OPF) problem is an extension of the
economic dispatch problem introduced by Carpentier in 1962
\cite{Car:62}. The goal of the OPF problem is to find a cost-optimal
operating-point of a power system that consists of a set of power
busses that are interconnected through a network of transmission
lines. Each power bus may have one or more generators and/or a load
(\ie, demand), and typical problem formulations minimize generation
cost or transmission loss, subject to a set of nonlinear power flow
constraints and bounds on \eg\ generation, voltage magnitudes, and
transmission line flows. Many different problem formulations exist,
but AC OPF problems are generally difficult nonconvex optimization
problems; see \eg\ \cite{Mom:01,Zhu:09} for further details on
different formulations and applications.

Numerous solution techniques have been proposed in previous work,
including nonlinear optimization techniques such as sequential
quadratic programming, Lagrangian relaxation, and interior-point
methods, and more recently, derivative-free methods such as particle
swarm optimization, genetic algorithms, and evolutionary programming;
see \eg\ \cite{QDB:09, ZMT:11} and references therein. Although the
derivative-free methods are generally more versatile, the conventional
nonlinear optimization techniques have some important advantages. For example,
derivative-free methods do not compute dual variables which have
valuable economic meanings in electricity markets.

In recent work, researchers have applied semidefinite relaxation (SDR)
techniques to various OPF problem formulations
\cite{BWFW:08,LaL:10,LaL:12}.  This approach leads to conic
optimization problems that can be solved in polynomial time with an
interior-point method. The SDR formulation has attracted significant
attention since, unlike previous methods, the solution provides either
(i) a certificate of global optimality, (ii) a certificate of
infeasibility, or (iii) a lower bound on the optimal value. We say
that the SDR is ``exact'' if it provides the global solution to the
original problem. Conditions that guarantee exactness of the SDR have
been derived and studied in \cite{BGLC:11,LaL:12,ZaT:13}.

The computational cost of solving the SDR problem grows rapidly with
the size of the OPF problem, and this renders the direct SDR
formulation impractical for large-scale OPF problems. This can be
attributed to two computational bottlenecks: (i) a dense symmetric
matrix variable that grows with the number of busses in the power
network, and (ii) a large dense positive definite system of equations,
the so-called \emph{Schur complement system}, that defines the search
direction at each iteration. Despite its apparent high computational
complexity, the SDR formulation possesses both sparsity and low-rank
structure.  The sparsity is related to the fact that each power bus is
typically connected to only a small number of adjacent power busses,
and, in fact, the dual variable inherits its sparsity pattern from the
network graph. This was pointed out by Lavaei \& Low \cite{LaL:12} who
proposed to solve the dual problem; however, they solve the dual
problem with SeDuMi \cite{Stu:99}, which is a general-purpose
primal--dual conic solver. The solver therefore maintains and factors
at each iteration not only the sparse dual variable, but also the
dense primal variable. It is important to note that the second
computational bottleneck (the dense Schur complement system) remains
even if one avoids forming the dense primal variable, for example, by
solving the SDR problem and/or its dual using a sparse semidefinite
programming (SDP) solver such as DSDP \cite{BeY:04} or SMCP
\cite{ADV:10}.

Inspired by advances in sparse semidefinite programming, Jabr
\cite{Jab:12} applied the conversion method by Fukuda \etal\
\cite{FKMN:00} to the SDR of the OPF problem. The conversion method is
based on positive semidefinite matrix completion, and it reformulates
a sparse SDP as an equivalent block-diagonal SDP with additional
equality constraints. When applied to a tree network, the conversion
method yields a block-diagonal SDP with $2\times 2$ blocks. In
general, the blocks correspond to the so-called cliques of a chordal
embedding of the network graph, and the extra equality constraints
impose consistency by forcing certain elements of different blocks to
be equal. Although this method typically introduces sparsity in the
Schur complement system, it also increases the order of the Schur
complement system quite substantially in many cases. The benefit of
using the conversion method therefore depends very much on the network
graph and on the chordal embedding. In \cite{MHLD:13}, Molzahn \etal\
apply the conversion method in combination with a simple heuristic for
reducing the number of cliques in the converted problem, and their
method typically reduces the computation time by a factor of
three. However, like previous work in this area, the conversion is
performed based on the SDR of a real-valued formulation of the
problem, and as we explain in Section \ref{s-new-sdr}, this approach
adds more than twice as many equality constraints to the converted
problem than necessary (unless the special structure of the
real-valued problem is taken into account).

\paragraph{Contributions} The main goal of this paper is to propose
new formulations and relaxations of the AC OPF problem. Specifically,
we propose a class of new computationally cheaper SDRs of the OPF
problem which are obtained by dropping some of the equality
constraints introduced when converting the original SDR to a
block-diagonal one.  We also propose a primal formulation of the AC
OPF problem in which we model line flow constraints using second-order
cones (SOCs) of dimension three instead of positive semidefinite cones
of order three. This reduces the number of variables. Furthermore, by using
complex-valued voltage variables in our model instead of their
real-valued real and imaginary parts, we obtain an SDR that involves a
Hermitian matrix variable of order equal to the number of power busses
instead of a symmetric matrix variable of twice the order. By applying
the conversion method of Fukuda \etal\ \cite{FKMN:00} to this SDR, we
avoid introducing unnecessary equality constraints.

\paragraph{Notation} 
We denote by $\reals^n$ and $\complex^n$ the sets of real and complex
$n$-dimensional vectors, and $\reals^{m \times n}$ and $\complex^{m
  \times n}$ denote the sets of real and complex $m \times n$
matrices. The set $\mathcal H^n$ is the set of Hermitian matrices of
order $n$, and $\mathcal H_+^n$ denotes the set of positive
semidefinite matrices in $\mathcal H^n$. The matrix inequality $A
\succeq B$ means that $A-B$ is positive semidefinite, \ie, the
eigenvalues of $A-B$ are nonnegative. We denote by $\iu = \sqrt{-1}$
is the imaginary unit, $a^*$ denotes the complex conjugate of $a \in
\complex$, and $A^H$ denotes the Hermitian transpose of a matrix $A
\in \complex^{m \times n}$. Finally, $\Tr(A)$ denotes the trace of a square
matrix $A$, and $\Tr(B^HA)$ denotes the inner product of $A$ and $B$.

\paragraph{Outline}
The rest of the paper is organized as follows. We first introduce the
power system model and our formulation of the OPF problem and its SDR
in Section \ref{s-prob-sdr}. We then review chordal conversion in
Section \ref{s-chordal-conversion} and introduce some new SDRs in
Section \ref{s-new-sdr}.  This is followed by numerical experiments in
Section \ref{s-numerical-experiments}, and we conclude the paper in
Section \ref{s-conclusions}.

\section{Problem Formulation and Semidefinite Relaxation}
\label{s-prob-sdr}
We start this section by describing the power system model and
a variant of the AC OPF problem. We then propose a reformulation
of the problem and derive the associated SDR.

\subsection{Model}
The power system model consists of a network of power busses.  We
denote the set of power busses (nodes) by $\mathcal N$ and the set of
transmission lines (edges) by $\mathcal L \subseteq \mathcal N \times
\mathcal N$, \ie, $(k,l) \in \mathcal L$ if there is a transmission
line from node $k$ to node $l$. Transmission lines may be nonsymmetric
and there may be more than one transmission line between a pair of
nodes, so we model the network graph as a directed graph.  We denote
the number of nodes by $|\mathcal N|$ and the number of transmission
lines by $|\mathcal L|$. For each $k \in \mathcal N$, we define
$\mathcal G_k$ as the set of generators at node $k$, and $\mathcal G =
\bigcup_{k \in \mathcal N} \mathcal G_k$ denotes the set of all
generators in the network. We allow $\mathcal G_k$ to be the empty set
which corresponds to a power bus without generators. We also define a
set $\Lflow \subseteq \mathcal L$ of 
transmission lines with flow constraints (\ie, $\Lflow = \mathcal L$ if all transmission lines have flow
constraints).

We associate with power bus $k$ a complex bus voltage $\V_k$ and a
complex current $\I_k$, and we define two vectors $\V =
(\V_1,\V_2,\ldots,\V_n)$ and $\I = (\I_1,\I_2,\ldots, \I_n)$. The
currents and voltages satisfy the equation $\I = Y \V$ where $Y \in
\complex^{n \times n }$ is a so-called bus admittance matrix which
inherits its sparsity from the network graph, thus $Y$ is
generally very sparse. We will henceforth assume that the network
graph is connected.  The bus admittance matrix follows from Kirchhoff's current
law, and it can be computed from the problem data; refer to \eg\
\cite{ZiM:11} for details on how to compute the bus admittance matrix.
We denote the complex power generated by generator $g$ by $\Sg_g =
\Pg_g + \iu\Qg_g$, and similarly, $\Sd_k = \Pd_k + \iu \Qd_k$ denotes
the known demand or load at bus $k$.

The OPF problem takes the following form
\begin{subequations} \label{e-opf}
\begin{align}
  \label{e-opf-objective}
   \mbox{minimize} \quad \sum_{g \in \mathcal G} f_g(\Pg_g)
\end{align}
subject to $\I = Y \V$ and the constraints
\begin{gather}
\label{e-opf-pbg}
 \I_k^*\V_k = \sum_{g\in \mathcal G_k} \Sg_g - \Sd_k, \quad k \in \mathcal N\\
\label{e-opf-gen-p}
P^{\min}_g \leq \Pg_g \leq P^{\max}_g, \quad g \in \mathcal
  G \\
\label{e-opf-gen-q}
 Q^{\min}_g \leq \Qg_g \leq Q^{\max}_g, \quad g \in \mathcal
  G \\ 
\label{e-opf-vm}
 V_k^{\min} \leq |\V_k| \leq V_k^{\max}, \quad k \in \mathcal N \\
\label{e-opf-fc1}
| S^{\mathrm{fl}}_{k,l}(\V) | \leq S_{k,l}^{\max}, \quad (k,l) \in
\Lflow \\
\label{e-opf-fc2}
| S^{\mathrm{fl}}_{l,k}(\V) | \leq S_{k,l}^{\max}, \quad (k,l) \in \Lflow.
\end{gather}
\end{subequations}
The constraints \eqref{e-opf-pbg} are the power balance equations,
\eqref{e-opf-gen-p}-\eqref{e-opf-gen-q} are real and reactive power
generation constraints, \eqref{e-opf-vm} are voltage magnitude
constraints, and \eqref{e-opf-fc1}-\eqref{e-opf-fc2} are constraints
on transmission line flows. The variables are $\I,\V$, and $\Sg_g =
\Pg_g + \iu \Qg_g$ for $g \in \mathcal G$, and
$S^{\mathrm{fl}}_{k,l}(\V)$ denotes the complex power flow from bus
$k$ to bus $l$ which is a quadratic function of $\V_k$ and $\V_l$. We
will use the notation $S^{\mathrm{fl}}_{k,l}(\V) = \V^H T_{k,l}\V +
\iu \V^H \Tb_{k,l} \V$ where $T_{k,l}$ and $\Tb_{k,l}$ are Hermitian
and given (see \cite{ZiM:11} for details regarding the transmission
line model). The scalars $P_g^{\min}$, $P_g^{\max}$, $Q_g^{\min}$,
$Q_g^{\max}$, $V_k^{\min}$, $V_k^{\max}$, and $S_{k,l}^{\max}$ are
real and given, and the cost function $f_g : \reals \to \reals$
represents the fuel cost model of generator $g$ and can be any
so-called semidefinite representable convex function. This includes
the set of linear and convex quadratic representable functions; see
\eg\ \cite{BeN:01}. Here we will model the fuel cost curve as a convex
quadratic function of the form
\begin{align} \label{e-convex-quad-cost}
  f_g(\Pg_g) = \alpha_g \Pg_g^2 + \beta_g \Pg_g,
\end{align}
where the scalars $\alpha_g \geq 0$ and $\beta_g$ are given; see \eg\
\cite{Mom:01}. 

\subsection{Reformulation and Semidefinite Relaxation}
If we denote by $e_k$ the $k$th column of the identity matrix of order
$n$, we can express the left-hand side of each of the power balance equations
\eqref{e-opf-pbg} as $\I_k^*\V_k = \V^HY^He_ke_k^T\V$,
and hence we can eliminate the variables $\I_k$ and remove the equation $\I=Y\V$
from \eqref{e-opf}. The real power
generation inequalities \eqref{e-opf-gen-p} can be expressed as \[ \Pg_g
= P^{\min}_g + \Pg_g^l, \quad \Pg_g^l + \Pg_g^u =
P^{\max}_g-P^{\min}_g,\] where $\Pg_g^l$ and $\Pg_g^u$ are nonnegative
slack variables, and the equation $\Pg_g = P^{\min}_g + \Pg_g^l$ allows us
to eliminate the free variable $\Pg_g$. In a similar fashion, we may introduce
nonnegative slack variables $\Qg_g^l$ and $\Qg_g^u$ for each of the
inequalities \eqref{e-opf-gen-q}, and each of the voltage constraints
\eqref{e-opf-vm} can be expressed as
\[ V_k^{\min} + \nu_k^l = |\V_k|, \quad |\V_k| +\nu_k^u =
V_k^{\max},\] with nonnegative slack variables $\nu_k^l$ and
$\nu_k^u$. We now define two Hermitian matrices $\Yk =
(1/2)(Y^He_ke_k^T + e_ke_k^TY)$ and $\Ybk = -(\iu/2)(Y^He_ke_k^T -
e_ke_k^TY)$ so that the real and reactive power balance equations at
bus $k$ can be expressed as
\begin{align*}
  \sum_{g \in \mathcal G_k} (P_g^{\min} + \Pg_g^l) &= \V^H\Yk \V +
  \Pd_k, \\
  \sum_{g\in \mathcal G_k} (Q_g^{\min} + \Qg_g^l) &= \V^H \Ybk
  \V+\Qd_k.
\end{align*}
The convex quadratic cost functions \eqref{e-convex-quad-cost} can be
expressed using the epigraph formulation $f_g(\Pg_g) \leq t_g$ where
$t_g$ is an auxiliary variable, and hence the inequality can be
expressed as the quadratic constraint $ | \sqrt{\alpha_g} \Pg_g |^2
\leq t_g - \beta_g \Pg_g$. This, in turn, is equivalent to a SOC
constraint (see \eg\ \cite{BeN:01})
\[
  \|
  \begin{bmatrix}    
  1/2-t_g + \beta_g \Pg_g\\
  \sqrt{2\alpha_g} \Pg_g
  \end{bmatrix}\|_2 \leq 1/2+t_g - \beta_g\Pg_g.
\]

Using the above reformulations and slack variables, and if we partition $\mathcal G$ into a set of generators
$G^{\mathrm{lin}} = \{ g \in \mathcal G \,|\, \alpha_g = 0 \}$ with linear
cost and generators $G^{\mathrm{quad}} = \{ g \in \mathcal G
\,|\, \alpha_g > 0 \}$ with quadratic cost, the OPF problem \eqref{e-opf}
can be expressed as
\begin{subequations} \label{e-opf-rc}
\begin{gather}
    \mbox{minimize}  \quad  \sum_{g \in \mathcal G^{\mathrm{lin}}} \beta_g
    (P_g^{\min} + \Pg_g^l) + \sum_{g \in \mathcal
      G^{\mathrm{quad}}} t_g
\intertext{subject to}
\label{e-opf-rc-power-bal-p}
   \Tr(\Yk X) = \sum_{g \in \mathcal G_k} (P_g^{\min} + \Pg_g^l) -
   \Pd_k, \  k \in
   \mathcal N\\
   \Tr(\Ybk X) = \sum_{g \in \mathcal G_k} (Q_g^{\min} + \Qg_g^l) -
   \Qd_k, \  k
   \in \mathcal N  \\
   \Pg_g^l + \Pg_g^u = P^{\max}_g- P^{\min}_g, \  g \in \mathcal  G\\
   \Qg_g^l + \Qg_g^u = Q^{\max}_g- Q^{\min}_g, \  g \in \mathcal  G\\
   (V_k^{\min})^2 + \nu_k^l = X_{kk}, \ X_{kk} + \nu_k^u =
   (V_k^{\max})^2, \  k \in   \mathcal N \\
   \label{e-opf-rc-fc1}
   z_{k,l} =
   \begin{bmatrix}
     (S_{k,l}^{\max})  \\ \Tr(T_{k,l}X) \\ \Tr(\Tb_{k,l}X)     
   \end{bmatrix}, z_{l,k} = \begin{bmatrix}
     (S_{k,l}^{\max})  \\ \Tr(T_{l,k}X) \\ \Tr(\Tb_{l,k}X)     
   \end{bmatrix}, \  (k,l) \in \Lflow\\
\label{e-opf-rc-epi-qgen}
   w_g =
  \begin{bmatrix}
     1/2 + t_g - \beta_g (P_g^{\min} + \Pg_g^l)  \\
     1/2 - t_g + \beta_g (P_g^{\min} + \Pg_g^l)  \\
     \sqrt{2\alpha_g} (P_g^{\min} + \Pg_g^l)
   \end{bmatrix}, \  g \in \mathcal G^{\mathrm{quad}} \\
\label{e-opf-rc-lin1}
   \Pg_g^l, \Pg_g^u, \Qg_g^l,\Qg_g^u \geq 0, \  g \in \mathcal G \\
\label{e-opf-rc-lin2}
   \nu_k^l, \nu_k^u \geq 0, \  k \in \mathcal N\\
\label{e-opf-rc-soc2}
   z_{k,l} \in \mathcal K_q^3,\ z_{l,k} \in \mathcal K_q^3, \  (k,l)
   \in \Lflow \\
\label{e-opf-rc-soc1}
   w_g \in \mathcal K_q^3, \ g \in \mathcal G^{\mathrm{quad}} \\
   \label{e-opf-rc-rank1}
   X = \V\V^H.
\end{gather}
\end{subequations}
Here $\mathcal K_q^3 = \{ (t,x) \in \reals\times\reals^2 \,|\, t \geq
\|x\|\}$ denotes the SOC in $\reals^3$. 
The constraints \eqref{e-opf-rc-fc1} and \eqref{e-opf-rc-soc2}
correspond to the line flow constraints
\eqref{e-opf-fc1}-\eqref{e-opf-fc2}, and the constraints
\eqref{e-opf-rc-epi-qgen} and \eqref{e-opf-rc-soc1} correspond to the
epigraph formulation of the cost functions for generators with
quadratic cost. 

The only nonconvex constraint in \eqref{e-opf-rc} is the rank-1
constraint \eqref{e-opf-rc-rank1}, and we obtain an SDR of the problem
simply by replacing the nonconvex constraint $X = \V\V^H$ with the positive
semidefinite constraint $X \succeq 0$. The SDR problem is convex, and
its solution $X^\star$ provides a lower bound on the optimal value of
\eqref{e-opf}. Furthermore, if $X^\star$ has rank 1, we obtain a
globally optimal solution to \eqref{e-opf} by computing a rank-1
factorization $X^\star = \tilde \V \tilde \V^H$. Note that $X^\star$
only carries relative phase information since $\tilde \V \tilde \V^H
= \bar\V \bar \V^H$ for
any $\bar\V = \exp(\iu \theta)\tilde \V$ where $\theta \in [0,2\pi]$.

The SDR of \eqref{e-opf-rc} can be expressed as a ``cone linear
program'' (cone LP) of the form
\begin{align} \label{e-opf-conelp}
  \begin{array}{ll}
    \mbox{minimize}  
    & h^Tz \\
    \mbox{subject to} 
    & G^Tz + c = 0 \\
    & z \in \mathcal K
  \end{array} 
\end{align}
with variable $z$ and where $\mathcal K$ is the direct product of
cones. To see this, notice that the constraints
\eqref{e-opf-rc-power-bal-p}-\eqref{e-opf-rc-epi-qgen} are linear
equality constraints. There are a total of 
\begin{align*}
  r = 4|\mathcal N| + 2|\mathcal G| + 3(2|\Lflow| + |\mathcal G^{\mathrm{quad}}|)
\end{align*}
of these, and they correspond to the constraint $G^Tz + c
= 0$ in \eqref{e-opf-conelp}. Similarly, the cone constraint $z \in
\mathcal K$ in \eqref{e-opf-conelp} corresponds to the $4|\mathcal
G|+2|\mathcal N|$ nonnegativity constraints
\eqref{e-opf-rc-lin1}-\eqref{e-opf-rc-lin2}, the $2|\Lflow| +
|\mathcal G^{\mathrm{quad}}|$ SOC constraints
\eqref{e-opf-rc-soc2}-\eqref{e-opf-rc-soc1}, and the constraint $X
\succeq 0$ (\ie, the relaxed version of \eqref{e-opf-rc-rank1}). In other
words, $z$ represents all variables in the SDR of
\eqref{e-opf-rc}, and the cone $\mathcal K$ is given by
\[ \mathcal K = \reals_+^{4|\mathcal G|+2|\mathcal N|} \times \mathcal
K_q \times \mathcal H_+^{|\mathcal N|} \] where $\mathcal K_q =
\mathcal K_q^3 \times \cdots \times \mathcal K_q^3$ is the direct
product of $2|\Lflow| + |\mathcal G^{\mathrm{quad}}|$ SOCs.  Thus, the
total number of variables is equal to
\begin{align*}
 4|\mathcal G|+2|\mathcal N| + 3(2|\Lflow| +
|\mathcal G^{\mathrm{quad}}|) + |\mathcal{N}|^2
\end{align*}
where $|\mathcal{N}|^2$ is the number of scalar variables needed to represent
the Hermitian matrix variable $X$.

Before we present our new relaxations of the OPF problem, we note that
the problem in \eqref{e-opf-rc} is equivalent to a nonconvex quadratic
optimization problem, and hence it is also possible to derive
relaxations based on linear optimization and SOC
optimization; see \eg\ \cite{KiK:01} and references
therein. While these relaxations are computationally more tractable
than the SDR of \eqref{e-opf-rc}, they are also weaker in general,
which means that they may provide more conservative lower bounds than
the standard SDR. For example, it is possible to obtain a SOC
 relaxation of \eqref{e-opf-rc} simply by replacing the constraint
$X \succeq 0$ in the SDR of \eqref{e-opf-rc} by positive
semidefiniteness constraints on some or all $2\times 2$ principal minors of
$X$. These constraints can be expressed as SOC
constraints since if $W$ is Hermitian and of order 2, then
\[ W =
\begin{bmatrix}
  x & y^* \\
  y & w
\end{bmatrix} \succeq 0 \ \Leftrightarrow \ \|
\begin{bmatrix}
  x - w \\
  2\Re y\\ 
  2\Im y
\end{bmatrix} \|_2 \leq x+w.
\]
It is important to note that the positive semidefiniteness of the $2\times 2$ principal minors of
$X$ is only a necessary condition, but not a sufficient condition, for positive
semidefiniteness of $X$, so a SOC relaxation is
generally weaker than the SDR. 

\section{Chordal Conversion}
\label{s-chordal-conversion}
We begin this section with a quick review of chordal sparsity and
chordal cones; a more comprehensive treatment of these concepts can be found in \eg\ \cite{BlP:93,Gol:04,And:11}.

\subsection{Chordal Sparsity Patterns and Cones}
\begin{figure}
  \centering
  \begin{tikzpicture}[inner sep=0pt]
    \tikzstyle{mynode} = [draw,circle,inner sep=0pt,minimum size=5mm,font=\footnotesize]
    \begin{scope}\setlength{\arraycolsep}{3pt}
      \node (s) at (-2,0.5) {$\begin{bmatrix}
 \times & \times & \times & 0 \\[-0.4ex]
  \times & \times & \times& \times\\[-0.4ex]
  \times & \times & \times & \times \\[-0.4ex]
  0       &\times&\times &\times  
\end{bmatrix}$};
      \node[mynode] (1) at (0,1) {$1$};
      \node[mynode] (2) at (0,0) {$2$};
      \node[mynode] (3) at (1,1) {$3$};
      \node[mynode] (4) at (1,0) {$4$};
      \draw (1)--(2);\draw (2)--(3);\draw (3)--(4);\draw
      (1)--(3);\draw (2)--(4);
      \node[draw, rectangle,inner sep=3pt,rounded corners,minimum width=2cm,font=\footnotesize] (c1) at
      (3,1) {$\gamma_1 = \{1,2,3\}$};
      \node[draw, rectangle,inner sep=3pt,rounded corners,minimum width=2cm,font=\footnotesize] (c2) at (3,0) {$\gamma_2 = \{2,3,4\}$};
      \draw(c1)--(c2);
    \end{scope}
  \end{tikzpicture}
  \caption{Example of sparsity pattern, sparsity graph, and clique tree.}
  \label{fig:sparsity-graph}
\end{figure}
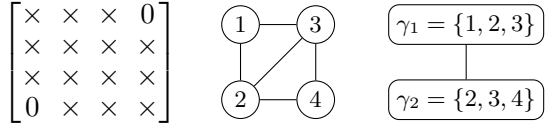

We represent a symmetric sparsity pattern of order $n$ as a set of
index pairs $E \subseteq \mathcal \{1,2,\ldots,n\} \times
\{1,2,\ldots,n \}$ where each pair $(i,j) \in E$ corresponds to a
nonzero entry of a sparse matrix of order $n$. We associate with the
sparsity pattern $E$ an undirected graph which has $n$ nodes
$\{1,2,\ldots,n\}$ and an edge between nodes $i$ and $j$ ($i\neq j$)
if $(i,j) \in E$. Fig.~\ref{fig:sparsity-graph} shows an example of a
sparsity pattern and the associated sparsity graph.
We say that the sparsity pattern $E$ is chordal if the sparsity graph
is chordal. A graph is chordal if all cycles of length greater than or
equal to 4 have a chord (\ie, an edge joining two non-adjacent
nodes in a cycle). Note that an immediate consequence of this
definition is that acyclic graphs are chordal. A clique of the
sparsity graph is a node subset $\gamma \subseteq \{1,2,\ldots,n\}$
such that the subgraph induced by $\gamma$ is complete and
maximal. Here ``maximal'' means that there exists no other complete
subgraph that properly contains the subgraph induced by $\gamma$. In
other words, each clique of the sparsity graph corresponds to a
maximal dense principal submatrix of the sparsity pattern $E$. It is
easy to verify that the sparsity graph in
Fig.~\ref{fig:sparsity-graph} is chordal (there are no cycles of
length 4 without a chord), and there are two cliques $\gamma_1 =
\{1,2,3\}$ and $\gamma_2 = \{2,3,4\}$. Indeed, these index sets
correspond to maximal dense principal submatrices in the sparsity
pattern. Note that the set $\{1,2\}$ also induces a complete subgraph,
but is is not maximal since it is contained in $\gamma_1$, and hence
it is not a clique according to our definition. As a final example, we
mention that a clique of a connected acyclic graph is any pair of
nodes that are connected by an edge, and this implies that a connected
acyclic graph with $n$ nodes has $n-1$ cliques of order 2.

We denote by $\mathcal E_\gamma (X)$ the
principal submatrix $X_{\gamma,\gamma}$ of $X$ defined by the index
set $\gamma$. Similarly, we denote the adjoint operator by $\mathcal
E_\gamma^{\mathrm{adj}} (W)$ which takes a matrix $W$ of order
$|\gamma|$ and returns a matrix of order $n$ with $W$ in the principal
submatrix defined by $\gamma$ and with zeros elsewhere.

Suppose $\gamma_1,\gamma_2,\ldots,\gamma_m$ are the cliques of a
connected chordal graph. (We have $m < n$ since a connected chordal
graph with $n$ vertices has at most $n-1$ cliques.) The cliques can be
arranged in a clique tree (a maximum weight spanning tree of the
clique intersection graph) that satisfies the so-called running
intersection property, \ie, $ \gamma_i \cap \gamma_j \subseteq
\gamma_k$ if clique $k$ is on the path between cliques $i$ and $j$ in
the tree; see \eg\ \cite{BlP:93}. Given a chordal graph, the cliques
and a clique tree can be found efficiently using \eg\ the algorithm by
Pothen \& Sun \cite{PoS:90}.  Nonchordal graphs can be handled by
means of a chordal embedding, and this technique is closely related to
sparse symbolic factorization techniques; see \eg\ \cite{ADV:12} and
references therein. In particular, a chordal embedding can be found
using a fill reducing reordering (such as ``approximate minimum
degree'' or ``nested dissection'') in combination with a symbolic
Cholesky factorization.

Given a sparsity pattern $E$, we define $\mathcal H_E^n$ as the set of
Hermitian matrices of order $n$ and with sparsity pattern $E$, and
$\mathcal P_E(X)$ denotes the projection of a (possibly dense)
Hermitian matrix $X$ onto $\mathcal H_E^n$. We define the cone of
positive semidefinite completable matrices in $\mathcal H_E^n$ as
$\mathcal H_{E,\mathrm{c}}^n = \{ \mathcal P_E(X) \,|\, X \succeq
0\}$.  A key result by Grone \etal~\cite[Theorem~7]{GJSW:84}
establishes that the cone $\mathcal H_{E,\mathrm{c}}^n$ is equivalent
to the set of partial positive semidefinite matrices in $\mathcal
H_E^n$ (\ie, matrices in $\mathcal H_E^n$ for which all dense
principal submatrices in $E$ are positive semidefinite) if and only if the
sparsity pattern $E$ is chordal.

\subsection{The Conversion Method}
The conversion method of Fukuda \etal\ \cite{FKMN:00} (which was
further studied and generalized by Kim \etal\ in \cite{KKMY:11}) makes
use of the aforementioned result by Grone \etal\ to express a cone
constraint $X \in \mathcal H_{E,\mathrm{c}}^n$ (where $E$ is
henceforth assumed to be chordal) as $m$ coupled constraints
\begin{align}\label{e-chordal-conversion-direct}
  W_k \succeq 0, \  W_k = \mathcal E_{\gamma_k}(X), \  k =
1,\ldots,m,
\end{align}
with $X \in \mathcal H_E^n$ and $W_k \in \mathcal H^{| \gamma_k|}$.
Using the running intersection property, we can eliminate the variable
$X$ in \eqref{e-chordal-conversion-direct}, \ie, $W_k = \mathcal
E_{\gamma_k}(X)$ for some $X \in \mathcal H_{E,\mathrm{c}}^n$ if and
only if for $k=1,\ldots,m$,
\begin{align}\label{e-chordal-conversion}
  W_k \succeq 0, \,  \mathcal E_{\gamma_j \cap \gamma_k}( \mathcal E_{\gamma_j}^{\mathrm{adj}}(W_j) -
  \mathcal E_{\gamma_k}^{\mathrm{adj}}(W_k)) = 0, 
  \,  j \in \mathrm{ch}(k)
\end{align}
where $\mathrm{ch}(k)$ is the set of indices of the cliques in the
clique tree that are children of clique $k$ (node $j$ is a child of
node $k$ in a rooted tree if there is an edge between $j$ and $k$, and
$k$ is on the path between $j$ and the root of the tree). Note that
the constraint $\mathcal E_{\gamma_j \cap \gamma_k}( \mathcal
E_{\gamma_j}^{\mathrm{adj}}(W_j) - \mathcal
E_{\gamma_k}^{\mathrm{adj}}(W_k)) = 0$ in \eqref{e-chordal-conversion}
couples principal submatrices of $W_j$ and $W_k$. These submatrices
are of order $|\eta_j|$ where $\eta_j = \gamma_j \cap \gamma_k$, and
since $W_k$ and $W_j$ are both Hermitian, the coupling consists of a
total of $|\eta_{j}|^2$ real equality constraints, \ie,
$|\eta_{j}|(|\eta_{j}|+1)/2$ equality constraints that correspond to the
symmetric real part and $|\eta_{j}|(|\eta_{j}|-1)/2$ equality
constraints that correspond to the skew-symmetric imaginary part. The
conversion \eqref{e-chordal-conversion} therefore introduces a total
of $s = \sum_{k=1}^m \sum_{j \in \mathrm{ch}(k)} |\eta_{j}|^2$
equality constraints. The Hermitian matrix inequality $X\succeq 0$ can
also be expressed as the real-valued symmetric matrix inequality
\[
Z = \begin{bmatrix}
  \Re X & - \Im X \\
  \Im X & \Re X
\end{bmatrix} \succeq 0
\]
and if we apply the conversion method to $Z$ without 
exploiting its particular structure, as has been done in previous work
(\eg, \cite{Jab:12,MHLD:13}), we add unnecessary consistency
constraints. Specifically, applying the conversion method directly to
$Z$ introduces $\sum_{k=1}^m \sum_{j\in \mathrm{ch}(k)} 2|\eta_{j}|
(2|\eta_{j}|+1)/2$ equality constraints, which is more than twice the
number necessary when conversion is applied to $X$.

The conversion method can be applied to a cone LP of the form
\eqref{e-opf-conelp}. However, since the conversion method only
affects the positive semidefinite matrix variables, we will simplify our
notation by considering the conversion method applied to a semidefinite
optimization problem of the form
\begin{align}
  \label{e-primal-csdp}
\begin{array}{ll}
  \mbox{minimize}  
  & \Tr(A_0X) \\
  \mbox{subject to} 
  & \Tr(A_jX) = b_j, \ j=1,\ldots,r, \\
  & X \in \mathcal H_{E,\mathrm{c}}^n
\end{array}   
\end{align}
where $A_j \in \mathcal H_{E}^n$, $j=0,1,\ldots,r$. Notice that we can express the inner
products $\Tr(A_jX)$ in terms of $W_1,\ldots,W_m$ as
\begin{align}
  \label{e-ai-splitting}
\Tr(A_jX) =
\sum_{k=1}^m \Tr(\widetilde A_{j,k} W_k).
\end{align}
where the matrices $\widetilde A_{j,k} \in \mathcal H^{|\gamma_k|}$ must
satisfy the condition $ A_j = \sum_{k=1}^m \mathcal E^{\mathrm{adj}}( \widetilde
A_{j,k})$. In general there are infinitely many ways of splitting
$A_j$ in this way. 

The conversion of \eqref{e-primal-csdp} is a semidefinite optimization
problem with $m$ blocks, \ie,
\begin{align}\label{e-primal-csdp-conversion}
  \begin{array}{ll}
    \mbox{minimize}  
    & \sum_{k=1}^m \Tr(\widetilde A_{0,k} W_k) \\
    \mbox{subject to} 
    & \sum_{k=1}^m \Tr(\widetilde A_{j,k} W_k) = b_j , \ j=1,\ldots,r,
    \\
    & \text{for }  k =1,\ldots,m \text{ and }  j \in \mathrm{ch}(k),\\
    & \quad \mathcal E_{\gamma_j \cap \gamma_k}( \mathcal E_{\gamma_j}^{\mathrm{adj}}(W_j) -
  \mathcal E_{\gamma_k}^{\mathrm{adj}}(W_k)) = 0 \\
    & W_k \succeq 0, \ k = 1,\ldots,m.
  \end{array} 
\end{align}
The total number of equality constraints is $r + s$ where $s$ is the
number of consistency constraints. For example, if the sparsity graph
is an acyclic graph, then there are $m = n-1$ cliques of order 2 and
$s = n-1$ consistency constraints. Thus, in this special case, the SDR
is equivalent to the SOC relaxation obtained by enforcing $n-1$
principal minors of $X$ to be positive semidefinite and expressing
these constraints as SOC constraints. For general problems, we will
see in Section \ref{s-numerical-experiments} that when applying the
conversion method directly to the cone LP \eqref{e-opf-conelp}, the
number of consistency constraints $s$ may be several times larger than
$r$. However, the number of variables in the converted problem is
typically much smaller than in the unconverted problem since the $m$
blocks are equivalent to $\sum_{i=1}^m |\gamma_i|^2$ scalar real-valued
variables instead of $|\mathcal N|^2$.

As mentioned in the introduction, the computational bottleneck when
solving an SDP is typically forming and factorizing the Schur
complement equations that define the search direction at each
interior-point iteration. The Schur complement matrix $H$ associated
with \eqref{e-primal-csdp} is of order $r$ with entries of the form
$H_{kl} = \Tr(A_kS^{-1}A_lS^{-1})$ for some $S \in \mathcal
H_{E,+}^n$, and hence $H$ is generally dense. The Schur complement
system associated with \eqref{e-primal-csdp-conversion}, however, is
often sparse, but it is of order $r+s$ instead of $r$. A detailed
exposition pertaining to conversion and sparsity in the Schur
complement system is outside the scope of this paper; see \eg\
\cite{KKK:08,SAV:13} and references therein.

\subsection{Clique Amalgamation}
The conversion method introduces a large number of equality
constraints when some of the sets $\eta_{j}$ are large. Amalgamating
(\ie, merging) a clique $j$ and its parent in the clique tree (say,
clique $k$), reduces the number of cliques by one, and it also reduces
the number of equality constraints by $|\eta_{j}|^2$. This, in turn,
implies that the order of the Schur complement system is reduced, but
it also affects the sparsity of the Schur complement system. Moreover,
the new combined clique is given by $\gamma_k \cup \gamma_j$, and
hence it is larger than both of the cliques from which it was
constructed. Thus, there are two different but coupled trade-offs to consider,
namely a trade-off between the number of cliques and their order (many
small cliques or fewer but larger cliques) and a trade-off between the
order of the Schur complement system and its sparsity (a large sparse
system or a smaller but less sparse system). 

In this paper, we will use the greedy clique amalgamation heuristic
from Sun \etal\ \cite{SAV:13} which does not take the sparsity of the
Schur complement system into account. Specifically, we start at the
bottom of the tree and merge clique $j$ and its parent clique $k$ if
$(|\gamma_k| - |\eta_{j}|)(|\gamma_j| - |\eta_{j}|) \leq
t_{\mathrm{fill}}$ or $\max(|\gamma_j|-|\eta_{j}|,
|\gamma_k|-|\eta_{k}|) \leq t_{\mathrm{size}}$. Here
$t_{\mathrm{fill}}$ is a threshold based on the amount of fill induced
by merging clique $j$ and its parent, and $t_{\mathrm{size}}$ is a
threshold based on the cardinality of the so-called supernodes which
are the sets $\gamma_k\setminus\eta_{k}$ and
$\gamma_{j}\setminus\eta_{j}$.

We end this section by noting that the term ``clique amalgamation'' is
inspired by terminology from the sparse factorization literature
where ``supernodal amalgamation'' refers to a similar technique used to balance the
number of supernodes and their orders to obtain cache-efficient block
sizes \cite{CR:89,HoSc:10}.

\section{Conversion-Based Semidefinite Relaxation}
\label{s-new-sdr}
In this section, we propose a new SDR technique for sparse nonconvex QPs
of the form
\begin{align}\label{e-nonconvex-qp}
  \begin{array}{ll}
    \mbox{minimize}  
    & x^H A_0 x \\
    \mbox{subject to} 
    & x^H A_j x = b_j, \ j=1,\ldots,r,
  \end{array} 
\end{align}
where $A_0,A_1,\ldots,A_r \in \mathcal H_E^n$ are given and $E$ is
chordal with cliques $\gamma_1,\ldots,\gamma_m$. An SDR of
\eqref{e-nonconvex-qp} is given by \eqref{e-primal-csdp} where $X \in
\mathcal H_{E,\mathrm{c}}^n$ is a convex relaxation of the constraint
$X=xx^H$. Recall that applying the conversion method to
\eqref{e-primal-csdp} yields the equivalent problem
\eqref{e-primal-csdp-conversion}. We refer to
\eqref{e-primal-csdp-conversion} as ``full conversion'' which is
closely related to the sparse SDR technique of Waki \etal\ \cite{WKKMM:06}
for polynomial optimization with structured sparsity.

Now recall that the consistency constraint
\eqref{e-chordal-conversion} that couples clique $j$ and its parent
clique is equivalent to $|\eta_{j}|^2$ equality constraints. By
dropping some of these equality constraints, we can obtain
computationally cheaper relaxations. We will refer to this method as
``conversion-based SDR'' (CSDR), and we will see in Section
\ref{s-numerical-experiments} that $O(|\eta_{j}|)$ equality
constraints linking each clique $j$ and its parent are often
sufficient to recover the solution to the standard SDR. This is not
surprising since we can verify whether two Hermitian rank-1 matrices
of order $|\eta_j|$ are equal or not by comparing a single column or
row. This corresponds to $2|\eta_{j}| -1$ equality constraints since
each row/column contains $|\eta_{j}|$ real values and $|\eta_{j}| -1$
complex values. Note that the CSDR solution is also a solution to the
standard SDR whenever it satisfies the full set of consistency constraints
\eqref{e-chordal-conversion}.
In the remainder of this section, we discuss CSDRs
based on different heuristic consistency strategies.

\subsection{Band CSDR}
A straightforward method for reducing the number of consistency
equalities is to keep only equalities that correspond to
entries of $\mathcal E_{\gamma_j \cap \gamma_k}( \mathcal
E_{\gamma_j}^{\mathrm{adj}}(W_j) - \mathcal
E_{\gamma_k}^{\mathrm{adj}}(W_k)) = 0$ within a band of half-bandwidth
$\rho_{j}$ (where $0 \leq \rho_{j} < |\eta_{j}|$). Diagonal consistency
corresponds to $\rho_{j} = 0$, and full consistency corresponds to
$\rho_{j} = |\eta_{j}|-1$. To simplify the choice of parameters, we
define a ``global'' half-bandwidth parameter $\rho \geq 0$ and let
$\rho_{j} = \min(\rho, |\eta_{j}|-1)$. The total number of consistency
constraints in the band CSDR is then given by
\[ \sum_{k=1}^m \sum_{j \in \mathrm{ch}(k)} \Bigl( |\eta_{j}| +
  2\sum_{l=1}^{\rho_{j}} (|\eta_{j}| - l) \Bigr).\] Recall that $\rho
= 0$ corresponds to diagonal-only coupling, and this is generally not
sufficient to recover a feasible solution to the original problem even
if the CSDR solution has rank-1 blocks. This is because the diagonal
elements of a rank-1 block correspond to squared voltage magnitudes,
and hence phase consistency is not enforced. However, for $\rho \geq
1$, we can always recover a solution to the original problem if all
blocks have rank 1 since the additional equality constraints are
sufficient to recover relative phase information.

\subsection{Other CSDRs}
An alternative to band CSDR is to keep only equality constraints that
correspond to diagonal entries as well as the entries in $\rho_{j} =
\min(\rho, |\eta_{j}|-1)$ rows/columns with $\rho \geq 0$ (\eg, an
arrow pattern). This approach leads to exactly the same number of
equality constraints as in the band CSDR with parameter $\rho$, and we
can recover a solution to the original problem if $\rho \geq 1$ and
all blocks have rank 1.

Yet another possibility is to keep equality constraints that
correspond to nonzero entries in the bus admittance matrix (\ie,
corresponding to edges in the network graph). The total number of
consistency equalities is then equal to $\sum_{k=1}^m \sum_{j \in
  \mathrm{ch}(k)} (|\eta_{j}| + 2|\mathcal L_j|)$ where $\mathcal L_j
\subseteq \mathcal L$
is the subset of transmission lines that connect a
pair of power busses that each belong to both clique $j$ and its
parent clique. We refer to this strategy as sparse CSDR.  It is also possible
to combine several consistency strategies, for example, using
band-plus-sparse structure.

As a final remark, we mention the possibility to use the conversion
technique based on the cliques of a nonchordal embedding of the
network graph (or the network graph itself) instead of the cliques of
a chordal embedding. Positive semidefiniteness of the cliques of a
partial Hermitian matrix is a necessary condition for it to have a
positive semidefinite completion. However, we know from the theorem of
Grone \etal\ \cite[Theorem~7]{GJSW:84} that this is not a sufficient
condition in general unless the sparsity pattern is chordal. Hence by
applying the conversion method based on the cliques of a nonchordal
patterns, we can obtain SDRs that are computationally cheaper to solve
(but generally also weaker) than the standard SDR.  We will not
explore this strategy further in this paper.

\section{Numerical Experiments}
\label{s-numerical-experiments}

We have implemented and tested the SDR of the problem in
\eqref{e-opf-rc} as well as some of the conversion-based SDR
techniques from Section \ref{s-new-sdr}. The experiments are based on
the benchmark problems from the \textsc{Matpower} package
\cite{ZiM:11}, and we build the data matrices associated with the cone
LP \eqref{e-opf-conelp} directly without using a modeling tool.  The
experiments were carried out in Matlab R2013a on a laptop with an
Intel Core i5 dual-core 1.8 GHz CPU and 8 GB RAM, and we used SeDuMi
1.3 with tolerance $\epsilon = 10^{-7}$ to solve the cone LP
\eqref{e-opf-conelp} and its different conversions. Note that although
we explicitly build the complex-valued cone LP \eqref{e-opf-conelp}
and apply the conversion method to this formulation, we 
cast the converted problem as a real-valued cone LP before passing
it to SeDuMi.

Following the approach in \cite{MHLD:13}, we treat generators with
tight upper and lower bounds as generators with fixed
power. Specifically, we set $\Pg_g = (P_g^{\min} + P_g^{\max})/2$ if
$P_g^{\max}-P_g^{\min}$ is less than 0.001 units for generator $g$,
and $\mathcal G^{\mathrm{fix}}$ denotes the set of number of such
generators. We also eliminate transmission line flow constraints that
are not active (\ie, operating below their maximum capacity) at the
(local) solution provided by \textsc{matpower}. (For the test problems
with several thousands of power busses, SeDuMi did not return a useful
solution when all the transmission line constraints were included in
the problem.)  Table \ref{tbl-benchmark-problems} lists the test cases
along with relevant problem dimensions. The number of active
transmission line constraints are listed in the last column of the
table.
\begin{table}
\center 
\begin{small} 
    \begin{tabular}{lrrrrrr}
      \hline
      Case  &         $|\mathcal N|$ & $|\mathcal G^{\mathrm{fix}}|$ &
      $|\mathcal G^{\mathrm{lin}}|$ & $|\mathcal G^{\mathrm{quad}}|$ &
      $|\mathcal L|$ & $|\Lflow|$\\ 
     \hline
      {IEEE-118} & 118  &  0   &    0 & 54 & 186 & 0 \\
      {IEEE-300} & 300  &  0   &    0 & 69 & 409 & 0 \\
      {2383wp}  & 2,383 & 92 & 235 & 0 & 2,896 & 5 \\
      {2736sp}   & 2,736 & 118 & 82 & 0 & 3,269 & 1 \\
      {2737sp}  & 2,737 & 165 & 54 & 0 & 3,269 & 1 \\
      {2746wop} & 2,746 & 346 & 85 & 0 & 3,307 & 0 \\
      {2746wp}  & 2,746 & 352 &104 & 0 & 3,279 & 0 \\
      {3012wp}  & 3,012 & 9 & 376 & 0 & 3,572 & 5 \\
      {3120sp}   & 3,120 & 25 & 273 & 0 & 3,693 & 8 \\
      \hline
    \end{tabular}
  \end{small}
\smallskip
\caption{Test cases and problem dimensions}
\label{tbl-benchmark-problems}
\end{table}

In order to improve the conditioning of the problem, we scale the
problem data $c$, $G$, and $h$ in the cone LP formulation
\eqref{e-opf-conelp} as $G := GD^{-1}$, $c := D^{-1}c$, and $h :=
h/\|h\|_2$ where $D= \diag(d_1,d_2,\ldots,d_{r+s})$
and $d_k = \max(|c_k|, \|Ge_k\|_{\infty} )$. This scales the dual
variables and the objective, but not the primal
variables. We found that this scaling heuristic reduced the number of
iterations and improved the solution accuracy in most of the
large-scale problems. Like in \cite{LaL:12}, we also use a minimum
line resistance of $10^{-4}$ per unit in our experiments.

For each test case, we solve the following SDRs: SDR with full
conversion based on chordal embedding with/without clique
amalgamation, band CSDR (with clique amalgamation), and sparse CSDR
(also with clique amalgamation). We use the clique amalgamation
parameters $t_{\mathrm{size}} = t_{\mathrm{fill}} = 16$.  Table
\ref{tbl-constraints} lists the number of constraints in each of the
SDR problems. The second column lists the number of constraints $r$ in
the unconverted problem, and columns 3--8 list $(r+s)/r$ which is the
total number of constraints \emph{normalized} by the number of constraints
in the unconverted problem.
\begin{table}
  \centering
    \begin{small}\setlength{\tabcolsep}{3pt}
    \begin{tabular}{lrcccccc}
      \hline
      Case & \# constraints & Full & Amal. & \multicolumn{3}{c}{Band
        CSDR ($\rho$)} & Sparse \\
      \cline{5-7}
            & (unconv.)& conv. &  conv.  & \multicolumn{1}{c}{$1$} &
              \multicolumn{1}{c}{$2$} & \multicolumn{1}{c}{$3$} & CSDR\\
      \hline
      {IEEE-118} & 742 $\times 1.0$  & 1.78 & 1.10 & 1.03 & 1.07 & 1.09 & 1.04  \\
      {IEEE-300} & 1,545 $\times 1.0$  & 1.99 & 1.14 & 1.04 & 1.09 & 1.12 & 1.04  \\
      {2383wp}   & 10,000 $\times 1.0$  & 2.97 & 1.63 & 1.08 & 1.20 & 1.30 & 1.10  \\
      {2736sp}   & 11,248 $\times 1.0$  & 3.04 & 1.64 & 1.08 & 1.22 & 1.32 & 1.10  \\
      {2737sp}   & 11,163 $\times 1.0$  & 3.03 & 1.62 & 1.08 & 1.22 & 1.32 & 1.10  \\
      {2746wop} & 11,379 $\times 1.0$  & 3.10 & 1.63 & 1.08 & 1.21 & 1.32 & 1.10  \\
      {2746wp}  & 11,438 $\times 1.0$  & 3.00 & 1.62 & 1.08 & 1.21 & 1.32 & 1.09  \\
      {3012wp}  & 12,716 $\times 1.0$  & 3.03 & 1.67 & 1.08 & 1.22 & 1.32 & 1.09  \\
      {3120sp}   & 12,990 $\times 1.0$  & 3.10 & 1.72 & 1.08 & 1.22 & 1.33 & 1.09  \\
      \hline
    \end{tabular}
  \end{small}
\smallskip
  \caption{Number of constraints before and after conversion}
  \label{tbl-constraints}
\end{table}
Recall that the ``full'' and ``amalgamated'' conversions are
equivalent to the unconverted problem whereas the ``band'' and
``sparse'' CSDRs are weaker, but computationally cheaper, relaxations.
The full conversion method adds the most equality constraints, and for
the large problems, this approach roughly triples the number of
constraints. The amalgamated conversion method is clearly much more
economical in term of the number of added equality constraints, but
the number of constraints still grows with more than 60\% when
converting large problems. The CSDRs introduce the smallest number of
constraints, \ie, around 10\%-30\% for large problems, depending on
the value of $\rho$. As a result, and as we will see later in this
section, the CSDRs are often much cheaper to solve.

\begin{table}
  \centering
  \begin{small}
    \setlength{\tabcolsep}{4pt}
    \begin{tabular}{lcccccc}
      \hline
      Case  & \multicolumn{1}{c}{Full} &
      \multicolumn{1}{c}{Amal.} & \multicolumn{3}{c}{Band CSDR ($\rho$)} & Sparse \\
      \cline{4-6}
            & \multicolumn{1}{c}{conv.} &  \multicolumn{1}{c}{conv.}  & \multicolumn{1}{c}{$1$} &
              \multicolumn{1}{c}{$2$} & \multicolumn{1}{c}{$3$} & \multicolumn{1}{c}{CSDR}\\
      \hline
      {IEEE-118} & $1.5\mbox{\textsc{e}+}{6}$ & $1.1 \mbox{\textsc{e}+}{7}$ & $1.1 \mbox{\textsc{e}+}{1}$ & $2.0 \mbox{\textsc{e}+}{7}$ & $1.9 \mbox{\textsc{e}+}{7}$ & $1.1 \mbox{\textsc{e}+}{1} $ \\
      {IEEE-300} & $2.4 \mbox{\textsc{e}+}{2}$ & $1.2 \mbox{\textsc{e}+}{3}$ & $1.7 \mbox{\textsc{e}-}{1}$ & $1.2 \mbox{\textsc{e}+}{3}$ & $1.1 \mbox{\textsc{e}+}{3}$ & $1.7 \mbox{\textsc{e}-}{1}$  \\
      {2383wp}   & $4.0 \mbox{\textsc{e}+}{2}$ & $5.3 \mbox{\textsc{e}+}{2}$ & $3.7 \mbox{\textsc{e}-}{1}$ & $9.4 \mbox{\textsc{e}+}{1}$ & $3.8 \mbox{\textsc{e}+}{2}$ & $4.1 \mbox{\textsc{e}-}{1}$  \\
      {2736sp}   & $ 1.6 \mbox{\textsc{e}+}{5}$ & $5.9 \mbox{\textsc{e}+}{5}$ & $1.2 \mbox{\textsc{e}-}{1}$ & $4.0 \mbox{\textsc{e}+}{2}$ & $1.1 \mbox{\textsc{e}+}{5}$ & $1.2 \mbox{\textsc{e}-}{1}$  \\
      {2737sp}   & $1.8 \mbox{\textsc{e}+}{4}$ & $7.2 \mbox{\textsc{e}+}{4}$ & $1.2 \mbox{\textsc{e}-}{1}$ & $6.2 \mbox{\textsc{e}+}{2}$ & $2.8 \mbox{\textsc{e}+}{4}$ & $1.2 \mbox{\textsc{e}-}{1}$  \\
      {2746wop} & $2.0 \mbox{\textsc{e}+}{4}$ & $5.6 \mbox{\textsc{e}+}{4}$ & $1.5 \mbox{\textsc{e}-}{1}$ & $1.2 \mbox{\textsc{e}+}{2}$ & $3.2 \mbox{\textsc{e}+}{4}$ & $4.8 \mbox{\textsc{e}-}{1}$  \\
      {2746wp}  & $3.5 \mbox{\textsc{e}+}{5}$ & $3.9 \mbox{\textsc{e}+}{5}$ & $1.3 \mbox{\textsc{e}-}{1}$ & $3.9 \mbox{\textsc{e}+}{2}$ & $1.9 \mbox{\textsc{e}+}{5}$ & $1.3 \mbox{\textsc{e}-}{1}$  \\
      {3012wp}  & $8.5 \mbox{\textsc{e}+}{0}$ & $1.6 \mbox{\textsc{e}+}{2}$ & $2.4 \mbox{\textsc{e}-}{1}$ & $7.5 \mbox{\textsc{e}+}{1}$ & $2.2 \mbox{\textsc{e}+}{2}$ & $2.4 \mbox{\textsc{e}-}{1}$  \\
      {3120sp}   & $6.6 \mbox{\textsc{e}+}{1}$ & $1.5 \mbox{\textsc{e}+}{2}$ & $2.8 \mbox{\textsc{e}-}{1}$ & $6.9 \mbox{\textsc{e}+}{1}$ & $1.6 \mbox{\textsc{e}+}{2}$ & $2.8 \mbox{\textsc{e}-}{1}$  \\
      \hline
    \end{tabular}
  \end{small}
\smallskip
  \caption{Eigenvalue ratios}
  \label{tbl-eig-ratio}
\end{table}
The ratio between the largest and the second largest eigenvalue can be
used as an indicator for the numerical rank of the solution. Roughly
speaking, the solution has numerical rank 1 if the aforementioned
eigenvalue ratio is sufficiently large. Since the converted problems
have multiple blocks, we consider the smallest such ratio, and these
are listed in Table \ref{tbl-eig-ratio}.  It is interesting to note
that full conversion with clique amalgamation yields slightly better
results than without. In particular, for the problem {3012wp}, there
is an order of magnitude difference between the eigenvalue ratio for
the two methods. Furthermore, for the band CSDR, the eigenvalue ratio
improves when the half-bandwidth $\rho$ is increased, and for $\rho =
3$, the eigenvalue ratios are comparable to those obtained via full
conversion. The sparse CSDR yields results that are similar to the
band CSDR with $\rho = 1$, and despite the lack of a rank-1
solution, the nearest rank-1 approximation may still be useful as
initialization for a general nonlinear solver.

Recall that in general, the CSDRs are weaker than the standard SDR. This
implies that the objective value can be used as an indication of the
relaxation quality. Table \ref{tbl-normalized-obj} lists the objective
values obtained via the CSDRs, normalized by the objective
values obtained via the standard SDR.
\begin{table}
  \centering
\begin{small}
  \begin{tabular}{lcccccc}
    \hline
      Case  & \multicolumn{3}{c}{Band CSDR ($\rho$)} & Sparse \\
      \cline{2-4}
      & \multicolumn{1}{c}{$1$} &\multicolumn{1}{c}{$2$} & \multicolumn{1}{c}{$3$} & CSDR\\
      \hline
       {IEEE-118} & 0.999 & 1.000 & 1.000 & 0.999  \\
      {IEEE-300}   & 0.999 & 1.000 & 1.000 & 0.999  \\
      {2383wp}   & 0.990 & 0.998 & 1.000 & 0.990  \\
      {2736sp}   & 0.989 & 1.000 & 1.000 & 0.990  \\
      {2737sp}   & 0.980 & 1.000 & 1.000 & 0.979  \\
      {2746wop} & 0.978 & 0.996 & 1.000 & 0.978  \\
      {2746wp}  & 0.989 & 1.000 & 1.000 & 0.989  \\
      {3012wp}  & 0.985 & 0.994 & 0.998 & 0.985  \\
      {3120sp}   & 0.988 & 0.999 & 1.000 & 0.989  \\
\hline
\end{tabular}
\end{small}
\smallskip
  \caption{Normalized objective value for CSDRs}
  \label{tbl-normalized-obj}
\end{table}
A normalized objective value of 1 corresponds to a relaxation that is
as tight as the original SDR. Notice that the CSDRs all yield lower
bounds that are within a few percent of the objective value obtained
via the standard SDR, and in all but one case, the band CSDR with
$\rho = 3$ yields a solution of the same quality as that obtained via
the standard SDR.

Finally we compare the computational complexity in terms of
computation time required to solve each of the relaxations. Table
\ref{tbl-cputime} shows ``wall time'' in seconds as reported by SeDuMi
(the ``CPU time'' was roughly a factor of two larger for all problems).
\begin{table}
  \centering
  \begin{small}\setlength{\tabcolsep}{3pt}
    \begin{tabular}{lrrrrrrrr}
      \hline
      Case & \multicolumn{1}{c}{No} & \multicolumn{1}{c}{Full} & \multicolumn{1}{c}{Amal.} & \multicolumn{3}{c}{Band
        CSDR ($\rho$)} & \multicolumn{1}{c}{Sparse}  & \cite{MHLD:13}\\
      \cline{5-7}
            & \multicolumn{1}{c}{conv.}  &\multicolumn{1}{c}{conv.} & \multicolumn{1}{c}{conv.}  & \multicolumn{1}{c}{$1$} &
              \multicolumn{1}{c}{$2$} & \multicolumn{1}{c}{$3$} & \multicolumn{1}{c}{CSDR} & \\
      \hline
      {IEEE-118} &  5.4 &   7.3 & 2.2 & 1.2 & 1.5 & 2.1 & 1.2 &
      2.1 \\
      {IEEE-300} &    78 & 19 & 5.0 & 3.8 & 4.2 & 4.0 & 3.7
      & 5.7 \\
      {2383wp}   & -& 650 & 225 & 78 & 103 & 132 &
      85 & 730 \\
      {2736sp}   & -&484 & 145 & 56 & 72 & 102 & 74
      & 622\\
      {2737sp}   & -&716 & 200 & 57 & 107 & 133 & 93 &607 \\
      {2746wop} & -&439 & 138 & 51 & 104 & 86 & 65 &738 \\
      {2746wp}  & -&547 & 168 & 50 & 95 & 95 & 82 &752 \\
      {3012wp}  & -&575 & 201 & 57 & 87 & 86 & 78 &1197 \\
      {3120sp}   & -&718 & 217 & 57 & 85 & 96 & 96 &1619 \\
      \hline
    \end{tabular}
  \end{small}
\smallskip
  \caption{Computation time}
  \label{tbl-cputime}
\end{table}
If we compare the full conversion and the amalgamated conversion, we see
that clique amalgamation typically results in a speed-up of around
2.5-3.5. Furthermore, the time required to solve the band CSDR with
$\rho = 3$ is roughly 100 seconds for the large-scale problems which
is a speed-up of around 5-6 when compared to the
time required to solve the full conversion SDR. The last column in the
table lists the results reported in \cite{MHLD:13},
and although a direct comparison is difficult (the experiments were
performed on different machines and using different formulations), the
large margins suggest that our approach is competitive and often much
faster than previously proposed methods.

\section{Conclusions}
\label{s-conclusions}
We have proposed a new method for generating computationally cheaper
SDRs of AC OPF problems. The method is based on chordal conversion,
which, given a chordal embedding of the network graph, converts a
sparse semidefinite optimization problem of the form
\eqref{e-primal-csdp} into an equivalent block-diagonal problem
\eqref{e-primal-csdp-conversion} that includes a set of consistency
constraints. By including only a subset of the consistency constraints,
we can generate conversion-based SDRs with reduced computational cost,
but these SDRs may also be weaker than the standard SDR. The band CSDR
method from Section \ref{s-new-sdr} keeps only consistency equalities
associated with entries within a band of half-bandwidth $\rho$, and
our numerical experiments indicate that this strategy works
surprisingly well in practice. More specifically, the band CSDR with
$\rho = 3$ has the same objective value as the standard SDR in all but
one test case, and the weaker band CSDRs yield objective values
that are within a few percent of those obtained via the standard
SDR. However, the experiments are based on only a small number of test
cases, so further experiments are necessary to thoroughly evaluate the
quality of the CSDRs.

In addition to the complexity-reducing CSDR technique, we lower the
computational cost further by applying the conversion method to
the complex-valued problem instead of its real-valued
counterpart. Moreover, we model transmission line flow constraints and
generators with quadratic fuel cost using SOC
constraints, which reduces the total number of variables and therefore
also computational cost. By combining these techniques, we have shown
that it is possible to solve SDRs of large-scale OPF problems significantly
faster than with previously proposed methods.

The CSDR technique can also be applied to extensions of the OPF
formulation \eqref{e-opf-rc}. One such extension is the so-called
multi-period OPF where power generation and demand vary over time, and
the time slots are coupled because of generator ramp rate limits. With
this kind of formulation, the standard SDR has a Hermitian positive
semidefinite matrix $X^{(i)}$ of order $|\mathcal N|$ for each time
slot $i$, and hence the problem dimension grows quickly with the
number of time periods. The CSDR technique can therefore be expected
to provide significant computational savings when applied to the SDR
of the multi-period OPF problem.

\end{document}